\documentclass[11pt]{amsart}
\usepackage{geometry}                
\usepackage{amssymb}
\usepackage{epsfig}
\usepackage{graphics,color,epic}

\input{epsf.sty}

\newcommand\pil[1]{\overset{#1}{\longleftrightarrow}}
\newcommand\kortpil{{\leftrightarrow}}
\newcommand\hpil[1]{\overset{#1}{\longrightarrow}}
\newcommand\tillmed[1]{\mathrel{ _{\rightarrow}{#1}}}
\newcommand\tillmot[1]{\mathrel{ _{\leftarrow}{#1}}}
\def\cu{{\mathcal U}}
\def\bp{{\bold p}}

\newtheorem{theorem}{Theorem}[section]
\newtheorem{lemma}[theorem]{Lemma}
\newtheorem{proposition}[theorem]{Proposition}

\newtheorem{remark}[theorem]{Remark}

\newtheorem{conjecture}[theorem]{Conjecture}

\begin{document}

\title{On percolation and the bunkbed conjecture}

\author{Svante Linusson} \thanks{ Svante Linusson is a Royal Swedish Academy of Sciences 
Research Fellow supported by a grant from the Knut and Alice Wallenberg 
Foundation.
}
\address{Department of Mathematics, KTH-Royal Institute of Technology, 
  SE-100 44, Stockholm, Sweden.}

\email{linusson@math.kth.se}

\date{November 30, 2009}

\begin{abstract}
We study a problem on edge percolation on product graphs $G\times K_2$. Here $G$ is any finite graph and $K_2$ consists of two vertices $\{0,1\}$ connected by an edge. Every edge in $G\times K_2$ is present with probability $p$ independent of other edges.  
The Bunkbed conjecture states that for all $G$ and $p$ the probability that $(u,0)$ is in the same component as $(v,0)$ is greater than or equal to the probability that 
$(u,0)$ is in the same component as $(v,1)$ for every pair of vertices $u,v\in G$.

We generalize this conjecture and formulate and prove similar statements for randomly directed graphs. The methods lead to a proof of the original conjecture for special classes of graphs $G$, in particular outerplanar graphs. 
\end{abstract}

\maketitle

\section{Introduction} \label{S:Intro}
This note is concerned with discussing a property of edge-percolation on finite graphs that should be intuitively clear, but more difficult to prove rigorously. To the best of my knowledge, the conjecture was first  formulated in a slightly different form (equivalent to model $E_2^p$ below) by P.W. Kasteleyn in 1985, see Remark 5 in \cite{vdBK}. In the form stated above, the conjecture has been presented in \cite {OH1} and \cite{OH2}, by Olle H\"aggstr\"om (who claimed it to be folklore). 

For any graph $G=(V,E)$ we consider the {\em bunkbed graph} $\tilde G:=G\times K_2$, where $K_2$ is the graph with two vertices $\{0,1\}$ and one edge. A vertex $x\in V(G)$ will have two images $x_0,x_1\in V(\tilde{G})$ and one edge between them. Such edges will be called {\em vertical} edges. Every edge $e\in E(G)$ has also two images $e_0,e_1\in E(\tilde{G})$ that will be called {\em horizontal} edges.
We will use the terms downstairs and upstairs to denote all vertices and edges in the 0-layer and 1-layer respectively. Let $0\le p\le 1$.

\medskip
\noindent\textbf{Model $E_1^p$} (Edge percolation):  Every edge in $\tilde{G}$ is present with probability $p$ independently of the other edges. We call the corresponding random graph $E_1^p(\tilde{G})$.

\smallskip
Of course this definition is not restricted to bunkbed graphs. For the theory of percolation in general we refer the reader to \cite{GG}.
For any bunkbed graph $\tilde{G}$, any vertices $x,y\in V(\tilde G)$ and $0\le p\le 1$, we define 
\begin{align*} P(x \pil {E_1^p(\tilde{G})} y):=& \text{ Probability that there is a path from $x$ to $y$ in 
$\tilde{G}$ under model $E_1^p$}.
\end{align*}
We will often omit $\tilde{G}$ if it is clear from the context what graph we are considering.
We will only be interested in connected graphs $G$.

The bunkbed conjecture $BBC_1^p(G)$ may now be defined as follows.

\begin{conjecture} [Bunkbed conjecture \cite{OH2}]
Let $G$ be any graph and $\tilde{G}=G\times K_2$ the corresponding bunkbed graph. For any $u,v\in V(G)$ and any $0\le p\le 1$ we have
\[ P(u_0 \pil{E_1^p(\tilde{G})} v_0)\ge P(u_0 \pil{E_1^p(\tilde{G})} v_1).
\]
\end{conjecture}

One motivation for formulating this problem is that we would like the probability $P(x\pil {E_1^p(\tilde{G})} y)$ to be a measure of how close the vertices $x$ and $y$ are in the graph $\tilde{G}$. For this to be a good concept we would like to make sure that "intuitive obvious" properties of closeness are true. To this end, the bunkbed graphs are natural testing candidates and we certainly want $BBC_1^p(G)$ to be true.
In \cite {OH2} H{\"a}ggstr{\"o}m coins the term bunkbed graph and proves the corresponding statement for a related model called random cluster model (also known as Fortuin-Kasteleyn model) with a certain parameter $q=2$. In that model graphs with a large number of non-connected components occur with higher probability and there is thus dependence between edges.
In \cite{BB}, Bollob\'as and Brightwell consider random walks on Bunkbed graphs and more general product graphs. They prove a number of interesting intuitively pleasing statements, but they also have a warning example that intuition sometimes might go wrong. In \cite{OH1} H{\"a}ggstr{\"o}m study continuous random walks on bunkbed graphs and proves a conjecture by Bollob\'as and Brightwell. The interesting papers \cite{vdBK} and \cite{vdBHK} have been inspired by the conjecture.

The rest of this paper is organized as follows. In Section \ref{S:Gen} we generalize the model in several steps to be able to use the combinatorial tools we want. In Section \ref{S:Outer} we prove (generalization of) $BBC_1^p(G)$ for outerplanar graphs. In Section \ref{S:Digraphs} we present the corresponding problem for randomly directed graphs. Lemma \ref{L:equal} states that probabilities of existence of paths in $E_1^{1/2}(G)$ are equal to existence of directed paths in the randomly directed case thus giving a direct connection to directed graphs. In Theorem \ref{T:D2} the corresponding bunkbed property is proved for a related model. Finally, in Section \ref{S:Critical} we define a critical probability for finite graphs.

\noindent
{\bf Acknowledgement:} I thank Olle H\"aggstr\"om for inspiring discussions on the conjecture.  I also thank J\"orgen Backelin, Madeleine Leander and a very helpful anonymous referee for comments on an earlier version of the manuscript.

\section{Generalizations and tools} \label{S:Gen}
We start with generalizing the model in three steps. The first step actually consists of two. First we condition on which vertical edges are present in $\tilde{G}$. Second we replace $p$ with a vector 
$\bold {p}=(p_e)_{e\in E(G)}, 0\le p_e\le 1$ for every $e$, giving a probability for each edge of $G$.
We will call such a vector $\bold{p}$ a {\em probability vector} on $G$. For given $\bold{p}$ and $T$ we define the following model.

\medskip
\noindent\textbf{Model $E_2^{\bold{p},T}$:} Vertical edges in $\tilde G$ are present exactly at positions in $T$. For each $e\in E(G)$, the horizontal edges $e_0, e_1$ in $\tilde G$ are present with probability $p_e$. All events that different edges are present are independent.

The vertices in $T$ will be called {\em transversal}. The natural generalization of the bunkbed conjecture, let us call it $BBC_2^{\bold {p},T}(G)$ seems also very likely to be true.

\begin{conjecture} [Kasteleyn, \cite{vdBK}]\label{C:E2}
Let $G$ be any graph and $\tilde{G}=G\times K_2$ the corresponding bunkbed graph. For any $u,v\in V(G)$, any $T\subset V(G)$ and any probability vector $\bp$ we have
\[ P(u_0 \pil {E_2^{\bold{p},T}(\tilde{G})} v_0)\ge P(u_0 \pil{E_2^{\bold{p},T}(\tilde{G})} v_1).
\]
\end{conjecture}

This conjecture is the original conjecture as formulated by P.W. Kasteleyn, see Remark 5 of \cite{vdBK}. In fact, that beautiful paper was inspired by Kasteleyn's conjecture.


\begin{proposition}
Given a graph $G$ and a probability $0\le p\le 1$, let $\bold {p}=(p_e)_{e\in E(G)}, p_e=p$ for every $e$. 
If  $BBC_2^{{\bold p},T}(G)$ is true for all $T\subseteq V(G)$, then $BBC_1^p(G)$ is true for the same graph $G$ and same $p$.
\end{proposition}

\begin{proof}
For any vertices $x,y\in \tilde G$ we have
\[P(x \pil {E_1^{p}(\tilde{G})} y)=\sum_{T\subset V(G)} P(x \pil {E_2^{\bold{p},T}(\tilde{G})} y)\cdot p^{|T|}\cdot 
(1-p)^{|V(G)\backslash T|}.
\]
The proposition follows.
\end{proof}

With this formulation one may start to prove the Bunkbed Conjecture for certain sets $T$ and we will present two easy examples. Recall that a set $C\subseteq V$ is called a {\em cutset} for $G$ if $G\setminus C$ is disconnected. If $x,y\in V$ are in different components of $G\setminus C$, then $C$ is said to {\em separate} $x$ and $y$. 

\begin{lemma} \label{L:cutset}
If $T\subseteq V(G)$ contains a cutset of $G$ separating $u$ from $v$, or if $u\in T$, or if $v\in T$ then
$P(u_0 \pil {E_2^{\bold{p},T}}v_0)=P(u_0 \pil {E_2^{\bold{p},T}}v_1)$, in particular $BBC_2^{{\bold p},T}(G)$ is true.
\end{lemma}

\begin{proof} This is easily proved with a mirror argument. Let $C\subseteq T$ be the cutset. Let $E_C\subseteq E$ be the edges in the same component as $v$ in $G\setminus C$ together with the edges with one endpoint in that component and the other in $C$. For any configuration of present edges $F\subseteq E(\tilde G)$, let $F_C:=
F\cap (E_C\times \{0,1\})$. Define $F_C'$ as the mirror image of $F_C$, i.e. an edge is present upstairs in $F_C'$ if and only if it is present downstairs in $F_C$ and vice versa. Define $F'=F\backslash F_C \cup F_C'$ and it is clear that the two configurations $F$ and $F'$ have the same probability. Also if there is a path from some vertex in $C$ to $v_0$ in $F$ then there is a path from the same vertex of $C$ to $v_1$ in $F'$ and vice versa. Since any path from $u_0$ to some vertex in $C\subseteq T$ may continue both upstairs and downstairs we receive a matching between cases with paths to $v_0$ and to $v_1$ respectively. The lemma follows.
\end{proof}

\begin{lemma} \label{L:fewvertices}
If $|T|=0,1$ then $BBC_2^{\bold{p},T}(G)$ is true for any graph $G$ and any probability vector $\bp$.
\end{lemma}
\begin{proof}
If $|T|=0$ then it is clear. Assume $T=\{x\}$. Given two complementary events $A_1,A_2$, it will suffice to prove $P(u_0\pil v_0\mid A_i)\ge P(u_0\pil v_1\mid A_i)$ for $i=1,2$. Here, we let $A_1$ be the event that there is no path from $u_0$ to $x_0$. Then the probability of a path to $v_1$ is zero, so the inequality follows. Let $A_2$ be the event that there exists a path from $u_0$ to $x_0$. The probability of a path from $x_1$ to $v_1$ upstairs is at most as large as a the probability of a path from $x_0$ to $v_0$ downstairs, because conditioning on the existence of a path from $u_0$ to $x_0$ downstairs can only affect the probability positively (Harris' inequality on increasing events \cite{H,GG}).  
\end{proof}

\medskip
To introduce our next generalization, let us fix any particular edge $e\in E(G)$. If $0<p_e<1$ there are four possibilities in model $E_2^{\bold{p},T}$ which we group into three cases as follows.

\begin{enumerate}
\item $e_0,e_1$ are both present
\item $e_0 $ is present and $e_1$ is absent, or\\ $e_0$ is absent and $e_1$ is present
\item $e_0,e_1$ are both absent
\end{enumerate}

The intuitive idea is to condition on which case we belong to and to use that case (1) can be thought of as contracting $e$ and (3) as removing $e$. This is made precise in Proposition \ref{P:E3impliesE2}. The remaining case (2) leads to defining the following new model for a given set $T\subseteq V(G)$.

\medskip
\noindent\textbf{Model $E_3^{T}$:} Vertical edges exist exactly at positions in $T$. Every horizontal edge upstairs in $\tilde{G}$ is present independently with probability $1/2$ and otherwise the corresponding edge exists downstairs. But no horizontal edge exists both upstairs and downstairs.

The natural generalization of the bunkbed conjecture is the following.

\begin{conjecture} [$BBC_3^{T}(G)$]
Let $G$ be any graph and $\tilde{G}=G\times K_2$ the corresponding bunkbed graph. For any $u,v\in V(G)$ and any $T\subseteq V(G)$ we have
\[ P(u_0 \pil {E_3^{T}(\tilde{G})} v_0)\ge P(u_0 \pil {E_3^{T}(\tilde{G})}v_1).
\]
\end{conjecture}

Recall that $G'$ is a {\em minor} of $G$ if it can be obtained by deleting and contracting edges of $G$.
For $e\in E(G)$ we use $G\backslash e$ and $G/e$ for the graph obtained when deleting and contracting the edge $e$. When we say minor in the proposition below we mean the usual notion in graph theory where multiple edges have been removed. However, it will later sometimes be convenient to allow multiple edges and then it will be explicitly stated.

\begin{proposition}\label{P:E3impliesE2}
If $BBC_3^{T'}(G')$ is true for any minor $G'$ of $G$ and all $T'\subseteq V(G')$, then $BBC_2^{{\bold p},T}(G)$ is true for any $\bp,T$ and thus also $BBC_1^p(G)$ for any $0\le p\le 1$.
\end{proposition}
\begin{proof}
Assume that $\bp, T$ and $G$ are given and that 
$BBC_3^{T'}(G')$ is true for any minor $G'$ of $G$ and all $T'\subseteq V(G')$
We will now condition on the edges of $G$ one at a time and prove the proposition by induction over the number of non-conditioned edges. For a given edge $e$, we will in case (1) contract $e$, in case (3) delete $e$ and in case (2) leave $e$ in the graph and remember that it now appears either upstairs or downstairs in the corresponding bunkbed graph. When we contract an edge we will in this proof allow the creation of multiple edges, but loops are irrelevant and may be deleted. Also when we contract an edge $xy$ we let the new vertex $v_{xy}$ be in $T$ if at least one of $x,y$ are in $T$. This way the probabilities for existence of paths will be preserved. Note that we have no assumption on $v\neq u$ in the bunkbed conjectures.

Let $F\subseteq E(G)$  and let $H$ be any graph where we have conditioned on the edges in $E(G)\backslash F$. So in $\tilde H$ for $e\in F$ we have that $e_0,e_1$ will occur independently with probability $p_e$. For an edge $e\in E(H)\backslash F$ exactly one of $e_0,e_1$ is present 
in $\tilde H$ each with probability $1/2$. 
The inductive hypothesis is that the corresponding bunkbed conjecture is true for any such graph $H$. 
With slight abuse of notation we will talk of such a graph $H$ also when we mean the entire model with 
probabilities for all possible configurations in $\tilde H$.

The base case, when $F=\emptyset$, is a graph $H$ as in Model $E_3^T$ with the difference that there might be multiple edges in $H$. If there are no multiple edges in $H$, then $H$ is a minor of $G$ and we are done. Assume $e,f\in E(H)$ are multiple (parallel) edges, then we consider the following complementary events: $A_1$ is the event that both $e_0$ and $f_0$ are present or that both $e_1$ and $f_1$ are present, $A_2$ is the event that $e_0,f_1$ or $e_1,f_0$ are present. If we condition on being in case $A_1$ then we let $H'=H\backslash f$. Since edge $f$ is irrelevant in this situation the conditional connection probabilities in $H$ are the same as the connection probabilities in $H'$. If we condition on $A_2$ then we contract $e$ and $f$ and call it $H''$ (possibly creating new multiple edges). Again, the conditional connection probabilities in $H$ are the same as the connection probabilities in $H''$. As in the proof of Lemma \ref{L:fewvertices} it suffices to prove the bunkbed inequality in the cases $A_1,A_2$ or equivalently for $H'$ and $H''$. Since $H'$ and $H''$ have strictly fewer edges we can perform another induction, this time over the number of edges, and it follows that they satisfy the bunkbed inequality.  As base case for this induction over $E(H)$ we have the graphs with no multiple edges.

For the inductive step, let $H$ be any graph obtained by conditioning on the edges in $E(G)\backslash F$ and $e\in F$. Let $H_1,H_3$ be the graphs obtained by contracting and deleting the edge $e$ respectively and note that connection probabilities $H_1$ (resp. $H_3$) are equal to the conditional connection probabilities in case (1) for $H$ (case (3) respectively).  Let also $H_{2}$ be the graph such that exactly one of $e_0$ and $e_1$ is present in $\tilde H_2$, which similarly correspond to the case (2). Thus, for any vertices $x,y\in V(\tilde H)$ we have that 

\[P(x \pil {\tilde{H}} y)=p_e^2\cdot P(x \pil {\tilde{H_1}} y)+2p_e(1-p_e)\cdot P(x \pil {\tilde{H_{2}}} y)+
(1-p_e)^2\cdot P(x \pil {\tilde{H_3}} y).
\]

For $H_3$ and $H_{2}$ the non-conditioned edges are $F\backslash {e}$, for $H_1$ they form a subset of $F\backslash {e}$ (edges parallel to $e$ become loops and thus removed). In any case we have by induction that all three graphs satisfy the bunkbed inequality.
It follows that the bunkbed conjecture is true also for $H$.
\end{proof}

It might seem more difficult to prove a conjecture not only for the graph $G$ but also for all its minors, but if the line of reasoning is to show that a minimal counterexample cannot exist then it is no more difficult.
Model $E_3$ has the great advantage that we no longer have the parameter $\bp$. 

\medskip

Another advantage is that we may reformulate it in terms of edge colorings of the original graph $G$ as follows.

\medskip
\noindent\textbf{Model $E_3^{T}$ reformulated:} Let $T\subseteq V(G)$. Every edge in $G$ is colored either red or blue with equal probability. A walk in $G$ may change color only at a vertex in $T$.
\medskip

Here we think of a blue edge as existing upstairs (blue as in heaven) and a red edge being downstairs. Arriving to $v_0$ or $v_1$ is the same as arriving to $v$ along a red and blue edge respectively. Recall that a {\em walk} in a graph is more general than a path since it is allowed to revisit a vertex. 
It is an elementary fact from graph theory that there exist a walk between two vertices if and only if there exists a path between the same vertices.
We need to use the term walk in this model since we could use a vertex both going along red edges and later along blue edges or vice versa. In the non-colored models the probability of existence of a path and of a walk is of course the same. We will from now on use mostly this second formulation
of $E_3^T$, but for notational convenience we use $\longrightarrow v_0$ for a walk entering $v$ along a red edge (if $v\in T$ blue edge also legal) and $\longrightarrow v_1$ for a walk entering $v$ along a blue edge (if $v\in T$ red edge also legal). 

Our proof in Section \ref{S:Outer} requires in fact yet one more level of generalization. In this next model we assume that some edges forming a connected subgraph are required to have the same color. To this end we think of the edges $E$ as partitioned into disjoint subsets $U_1,\dots ,U_k$, i.e. $\cup_{i} U_i=E$ and $U_i\cap U_j=\emptyset$ if $i\neq j$.  Let $\cu =\{U_1,\dots ,U_k\}$ be such a partition into connected subgraphs and $T\subseteq V(G)$.

\medskip
\noindent\textbf{Model $E_4^{T,\cu}$:} (Hypergraph)  All edges in a set $U_i$ are given the same color red or blue with equal probability independent of the other sets. A walk in $G$ may change color only at a vertex in $T$.\\

\smallskip
Note that model $E_3$ is the special case, where all sets $U_i$ contains one edge.  It is helpful to think of a set $U_i$ as a hyperedge having a color, which enables passage between any two of the vertices in the hyperedge. Thus model $E_4$ is a generalization to hypergraphs.

\begin{conjecture} [$BBC_4^{T,\cu}(G)$]
Let $G$ be any graph and $T\subseteq V$ and $\cu$ as in model $E_4$. For any $u,v\in V(G)$ we have
\[ P(u_0 \pil {E_4^{T,\cu}(G)} v_0)\ge P(u_0 \pil {E_4^{T,\cu}(G)} v_1).
\]
\end{conjecture}

It seems me that also this more general conjecture is likely to be true. 

\begin{remark}\label{R:olika}
It is worth noting however that one may not in general assume that two edges have different colors without violating the bunkbed condition. As an example, let $G$ be the path of length two from $u$ to $v$ with $x$ as middle vertex and let $T=\{x\}$. If we now assumed that the two edges $ux,xv$ have to have {\em different} colors then we would have 
\[ 0=P(u_0 \pil{} v_0)<P(u_0\pil{} v_1)=1/2,
\]
contrary to what we conjecture in the other models.
\end{remark}

\section{Outerplanar graphs} \label{S:Outer}

In this section we will prove the Bunkbed conjectures $BBC_1^T(G), BBC_2^{{\bold p},T}(G)$ and $BBC_3^T(G)$ for outerplanar graphs $G$. A connected planar graph is called {\bf outerplanar} if it is has a drawing such that every vertex lies on the boundary of the outer region.
This is equivalent to the graph not having $K_4$ or $K_{2,3}$ as minors.

\medskip
Our line of proof is to recursively prove that a minimal counterexample may not exist. To this end we will present a number of recursive operations. We will often need to work in model $E_4^{T,\cu}$. In each case we have a triple  $(G,T,\cu)$, a graph $G$, a set of transversal vertices $T\subseteq V(G)$ and a partition $\cu$ of $E(G)$. We say that the triple {\em reduces} to a set of triples 
$(G_i,T_i,\cu_i)$ if whenever $BBC_4^{T_i,\cu_i}(G_i)$ is true for all $i$ then also
 $BBC_4^{T,\cu}(G)$ is true. The operations below will be constructed by conditioning on mutually exclusive events and thus every probability for a walk in $(G,T,\cu)$ is a linear combination of the probability of the corresponding walks in $(G_i,T_i,\cu_i)$ which implies that $(G,T,\cu)$ reduces to $(G_i,T_i,\cu_i)$.
When we are interested in $BBC_3$ only we take $\cu$ to be the partition into singletons.

Whenever we contract an edge $xy$ we let the new vertex $v_{xy}$ be in $T_i$ if at least one of $x,y$ are in $T$. To avoid technicalities we will allow multigraphs and keep multiple edges that are formed after a contraction. Loops will always be removed. Note that we have not assumed that $u\neq v$ in general. This might be the case after a contraction and it is no problem.

\medskip
\noindent\textbf{T-operation:} If $x,y\in T$ and $xy\in E(G)$, then we contract the edge $xy$ to the graph $G_1:=G/xy$, with $T_1:=T\backslash \{x,y\}\cup \{v_{xy}\}$ and $\cu_1$ the partition $\cu$ restricted to the new edge set. 

\smallskip
Any walk can always run freely between the vertices $x$ and $y$ and assume any color when leaving $x,y$. Thus every probability $P(u_0\pil {}v_i)$ is preserved when contracting the edge $xy$. 
When this is the case we will call the graphs {\bf equivalent}.
Thus $(G,T, \cu)$ reduces to $(G_1,T_1,\cu_1)$.

\medskip
\noindent\textbf{V2-operation:} Assume $x\in V\setminus (T\cup \{u,v\})$ and $\deg(x)=2$. Let $y,z$ be the neighbors of $x$ and assume that at least one of $xy$ and $xz$ form a singleton set in $\cu$. Then we define two subgraphs of $G$ as follows. $G_1=G\setminus x$ and $G_2:=G/xy$ with $T_1=T_2=T$ and the natural restrictions on $\cu$.

\smallskip
If the edges $xy,xz$ have different colors we are in a situation equivalent to $G_1$. If $xy,xz$ have the same color, then we are in a situation equivalent to $G_2$. Thus $(G,T,\cu)$ reduces to $(G_i,T_i,\cu_i)$, $i=1,2$.

\medskip
\noindent\textbf{$\Delta$-operation:} Assume $x,y,z\in V$ are any vertices (possibly including $u$ or $v$) of the graph that form a triangle, i.e. $xy,xz,yz\in E(G)$. Assume further that no other edge is dependent on the color of $xy,xz$ or $yz$, i.e. each of them form a  singleton set $U_i$. Then we form the following four cases: $G_1:=G/xy$, $G_2:=G/xz$, $G_3:=G/yz$ and finally $G_4$ is the same graph as $G$, but we require $xy,xz$ and $yz$ to have the same color so they form on block $U=\{xy,xz,yz\}$ in the partition $\cu_4$. In this particular situation we do not want the graph $G_1$ to have double edges between $v_{xy}$ and $z$ so we remove one of them and similarly for $G_2,G_3$. Other multiple edges could have been created as usual.

\begin{figure}[htbp]
\begin{center}
\epsfig{file=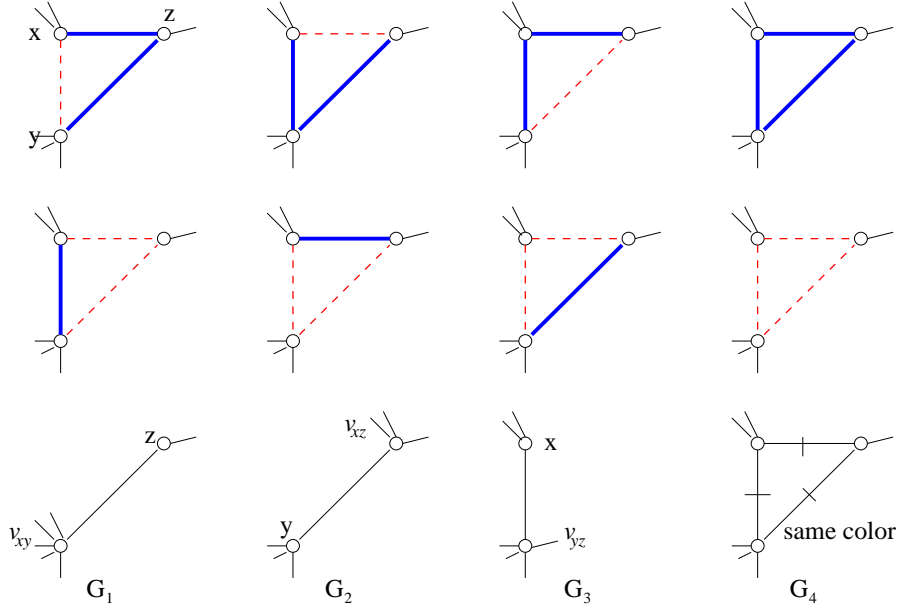, height=8cm}
\end{center}
\caption{The eight possible coloring of a triangle and how they are paired to reduce to simpler graphs $G_1,G_2,G_3$ and (hypergraph) $G_4$. Dashed edges are red, thick edges are blue.} 
\label{F:Delta}
\end{figure}

%

\smallskip
Assume first that $x,y,z\notin T$. Consider the eight possible colorings of $xy,xz,yz$, see Figure \ref{F:Delta}. Let case $A_1$ be the two leftmost figures, where $xy$ has a different color than $xz,yz$. Then we see that any walk can run freely between vertices $x$ and $y$ along either of the two different colors. This case is thus equivalent to $G_1$. Since we have assumed that $xz,yz$ have the same color they have this also in $G/xy$ and this is the reason we removed one of them in the definition of $G_1$. Similarly let $A_2$ be the case where $xz$ has a different color than $xy,yz$ and let $A_3$ be the case where $yz$ has a different color than $xy,xz$. Then these case are for the same reasons equivalent to $G_2$ and $G_3$ respectively. In the two rightmost figures the colors are equal for all three edges which is $G_4$. 
Thus $(G,T,\cu)$ reduces to $(G_i,T_i,\cu_i), i=1,2,3,4$. 

If any of $x,y,z$ belong to $T$, the reduction still works to the same four triples. For instance, if 
$x\in T$, $y,z\notin T$, let $A_1$ be the same case as above. If, say, $xy$ is red and $yz,xz$ are blue, then any walk entering $z$ blue can leave $x$ or $y$ in any color. Entering $x,y$ in either color the walk can leave at the other vertex in either color, or at $z$ in blue. Similarly with the colors reversed. Thus again, conditioning on case $A_1$ is equivalent to $G_1:=G/xy$. We leave to the reader to verify all other possibilities, which are not more difficult.

\begin{theorem}  $BBC_3^T(G)$ is true for all outerplanar graphs $G$ and all possible $T\subseteq V(G)$. Thus the bunkbed conjectures $BBC_1^p(G), BBC_2^{{\bold p},T}(G) $ are true for any outerplanar graph $G$ and any $p,\bp,T$.
\end{theorem}

\begin{proof} We will in fact prove the theorem for outerplanar multigraphs. 
Assume the contrary and let $G$ be a minimal counterexample, for some set $T$.
Minimal here means that all graphs obtained by deleting or contracting an edge are not counterexamples for any set $T$. Note that if $u\in T$ then we get equality by a mirror argument changing the color of every edge, similar to the proof of Lemma \ref{L:cutset}. Hence, we may assume that $u\notin T$. We may also assume that $G\backslash u$ is connected, since otherwise we could reduce to the component containing $v$. Similarly $G\backslash v$ is connected. In fact we may assume that $G$ is 2-connected. If not there would be a cutpoint $x$ such that $G\backslash x$ is disconnected. If there is a component $C$ s.t. $u,v\notin C$ then we can condition on the colors of the edges in $C$ which will imply a situation where we may or may not change color using a tour into $C$. This is equivalent to conditioning on if $x\in T$ or not which means that $G$ is not a minimal counterexample. If $u$ and $v$ are in different components $C_1,C_2$, then let $G_1:=G\backslash C_2$ and $G_2:=G\backslash C_1$. In this situation every walk from $u$ to $v$ passes through $x$ so 
\begin{align*}
&P_{G}(u_0\kortpil v_0)-P_{G}(u_0\kortpil v_1)=\\
&P_{G_1}(u_0\kortpil x_0)P_{G_2}(x_0\kortpil v_0)+
P_{G_1}(u_0\kortpil x_1)P_{G_2}(x_1\kortpil v_0)-P_{G_1}(u_0\kortpil x_0,x_1)P_{G_2}(x_0,x_1\kortpil v_0)-\\
\big(&P_{G_1}(u_0\kortpil x_0)P_{G_2}(x_0\kortpil v_1)+
P_{G_1}(u_0\kortpil x_1)P_{G_2}(x_1\kortpil v_1)-P_{G_1}(u_0\kortpil x_0,x_1)P_{G_2}(x_0,x_1\kortpil v_1)\big)=\\
&\big(P_{G_1}(u_0\kortpil x_0)-P_{G_1}(u_0\kortpil x_1)\big)
\big(P_{G_2}(x_0\kortpil v_0)-P_{G_2}(x_0\kortpil v_1)\big)\ge 0
\end{align*}
Here we use that by symmetry $P_{G_2}(x_0\kortpil v_1)=P_{G_2}(x_1\kortpil v_0)$, $P_{G_2}(x_0\kortpil v_0)=P_{G_2}(x_1\kortpil v_1)$ and 
$P_{G_2}(x_0,x_1\kortpil v_0)=P_{G_2}(x_0,x_1\kortpil v_1)$. The notation $P_{G_2}(x_0,x_1\kortpil v_1)$ means the probability that there are walks in $G_2$ from both $x_0$ and $x_1$ to $v_1$.

\medskip

So we may assume that $G$ is 2-connected and there are therefore two independent paths from $u$ to $v$ along the outer region, call them the outer paths. A {\em chord} is any edge not in the boundary of the outer region.

\noindent
\textbf{Claim 1:} All chords $xy$ in $G$ separates $u$ and $v$, that is $u$ and $v$ are in different components of the graph obtained by removing vertices $x,y$ from $G$.\\

Note that this implies in particular there are no chords with $u$ or $v$ as an endvertex.
To prove the claim we assume the opposite, that $G$ contains a chord between two vertices on the same outer path from $u$ to $v$. Then there is one such chord $xy$, with as few vertices $z_1,\dots,z_k$ as possible between $x$ and $y$ along the outer path. By construction $\deg(z_1)=\dots\deg (z_k)=2$ and $u,v\notin \{z_1,\dots, z_k\}$. If any $z_i\notin T$ then we can use operation $V2$ to reduce to smaller graphs for which the bunkbed conjecture is true by assumption, which gives a contradiction. Similarly if $z_i,z_{i+1}\in T$ we get a contradiction from the T-operation. This gives that the only possible configuration is a triangle $x,z,y$, where $z$ is of degree 2 between $x$ and $y$ along the outer path and $z\in T$. The $\Delta$-operation reduces to subgraphs $G_1,G_2,G_3$, for which the conjecture is true by assumption and $G_4$ where the three edges of the triangle have the same color. In the latter case one may remove $z$ and its two edges without altering any probability. This is again a subgraph of $G$ with no color assumptions and this contradicts $G$ being a minimal counterexample. The claim follows.

\medskip
The claim has the direct consequence $\deg(u)=2$. Let $x,y$ be the neighbors of $u$. We now condition on the color of $ux$.  If it was blue (corresponding to upstairs) we can never use that 
edge for any walk 
containing $u_0$ (downstairs) since $\deg(u)=2$. In that case we could remove $ux$ to obtain a smaller graph for 
which the theorem is true by assumption. We may thus assume that $ux$ is red and is a minimal 
counterexample when $\deg(u)\le 2$ and $ux$ is red. Similarly we can argue that $uy$ is red.

If $xy\notin E(G)$ then Claim 1 and the outerplanarity of $G$ implies that one of $x,y$ say $y$ has degree two and we may contract the red edge $uy$ to obtain a minor $G'$. This graph is a smaller graph with $\deg(u)\le 2$ and the condition that edge $ux$ is red, which by assumption is not a counterexample.

If $xy\in E(G)$, then we condition on the color of $xy$.  Again, because of outerplanarity and Claim 1 one of $x,y$, say $y$ has degree at most 3. If $xy$ is blue, then $x$ and $y$ are connected both with a red path and a blue edge. We may thus contract $xy$ without changing any probabilities for walks. Since $ux$ and $uy$ are both red no path can ever enter $u_1$ and every path starting in $u_0$ must first go to $v_{xy}$. We may thus contract also $ux,uy$ and the resulting minor must satisfy the bunkbed conjecture. If $xy$ is red, then we may contract $uy$ and remove one of the parallel red edges $uy,xy$ to get a new graph $G_1$. Probabilities for walks starting in $u$ in $G$ will be the same as walks starting in $u':=v_{uy}$ in $G_1$. Since $\deg(y)\le 3$ we get $\deg_{G_1}(u')\le 2$ and
$G_1$ has exactly one red edge $u'x$. But $G$ was a minimal such counterexample so the bunkbed conjecture is true for $G_1$ and we get the desired contradiction.
\end{proof}

\medskip
Note that there are other operations that one possibly may use to prove the conjectures for larger classes of graphs. We end this section with two examples.

\medskip
\noindent\textbf{Restricted $\Delta$-operation:}  Assume $x,y,z\in V$ form a triangle, i.e. $xy,xz,yz\in E(G)$. Assume further that $xy\in U_i$, $|U_i|\ge 2$, whereas the color of $xz$ and $yz$ is not dependent on the color of any other edge. Then we form the following three cases: $G_1:=G/xz$,  $G_2:=G/yz$ and finally $G_3$ is the same graph as $G$, but we require $xz$ and $yz$ to have the same color. As in the $\Delta$-operation we remove the multiple edge $yz$ in $G_1$ and the edge 
$xz$ from $G_2$. The set $U_i\in \cu$ such that $xy\in U_i$ do not change. 
The same reasoning as for $\Delta$-operation shows that $G$ reduces to $G_i,i=1,2,3$. 

The reason we cannot use the ordinary $\Delta$-operation is that if we contract $xy$ this would form a situation where the edges $U_i\backslash xy$ are forced to have different color than $yz$, which is not legal in model $E_4$. See also Remark \ref{R:olika}.

\medskip
\noindent\textbf{Y-operation:} Assume $x\in V\setminus T$ and $\deg(x)=3$. Let $a,b,c$ be the neighbors of $x$ and assume that the color of no other edge is dependent on the color of $ax,bx,cx$.
Then we form four subgraphs of $G$ as follows. $G_1=(G\setminus ax) /bx$, 
$G_2=(G\setminus bx) /cx$, $G_3=(G\setminus cx) /ax$ and $G_4$ is the same graph as $G$ but the edges $ax,bx,cx$ must have the same color.

\smallskip
If the edges $bx,cx$ have the same color but different from $ax$ we are in a situation equivalent to $G_1$. If the color of $bx$ is different from $ax,cx$ then we are in a situation equivalent to $G_2$ and similarly for $G_3$. The remaining cases are when all three edges have the same color which gives $G_4$. As for previous operations we see that $G$ reduces to $G_i,i=1,\dots,4$.\\ 
There is also a restricted Y-operation, whose formulation is left to the reader.

\medskip
Note that a unicolored $Y$ and a unicolored $\Delta$ give the same hypergraph. This opens the possibility to perform $\Delta \kortpil Y$ transformations of graphs. It is well-known that every planar graph is $\Delta\kortpil Y$ reducible to $K_2$. I have however not been able to use this fact to 
prove the $BBC_3^T(G)$ for planar graphs. One obstacle is that one may perform the Y-operation only if $x\notin T$.

\section{Randomly oriented graphs} \label{S:Digraphs}
In this section we present a connection to randomly directed graphs. First the basic model.

\medskip
\noindent\textbf{Model $D_1$:}  Every edge in ${G}$ is given one of the two possible directions with equal probability independently of the other edges.\\
We call the corresponding random directed graph $D_1({G})$.

By analogy with the undirected case we define $ P(x \hpil {D_1({G})} y):=$ Probability that there exist a directed path from vertex $x$ to $y$ in $G$ under model $D_1$. 
This model is a natural candidate to define a random orientation of a given graph. It was for example studied for the 
$\mathbb Z^2$-lattice in \cite{Gr2000} and for questions of correlation of directed paths in \cite{AL1,AL2}.

The following lemma gives a direct connection between model $D_1$ and $E_1^{1/2}$. It gives an interesting non-trivial reformulation of the problem. It is, to the best of my knowledge, first published by McDiarmid \cite{CM} and seemingly independently and with an elegant proof by Karp \cite{K} (My thanks to Jeff Kahn and the anonymous referee for pointing out these two references.) The lemma might seem surprising at first sight but once discovered it is not so difficult to prove. A third proof can be found in \cite{SL}.

\begin{lemma} \label{L:equal}
For any graph $G$ and any vertices $x,y\in V(G)$ we have
\[ P(x\pil {E_1^{1/2}(G)}y)=P(x\hpil {D_1(G)} y).
\]
\end{lemma}

This means that for the special case $p=1/2$, we may study randomly oriented graphs instead. A different model of directed graphs that also is applicable to other values of $p$ is discussed in \cite{SL}. Note that $D_1$ is a truly different model than $E_1$. We may for instance not generalize by conditioning on the direction of vertical edges as we have conditioned on the presence of vertical
 edges in $E_2$. We have however not been able to prove the bunkbed conjecture using these directed graphs either. 

The reformulation of model $E_3$ inspired the following two models replacing red and blue with directions.

\medskip
\noindent\textbf{Model $D_2^{T}$:} Let $T\subseteq V$ . Every edge in $G$ is given one of the two possible directions with equal probability. A walk in $G$ may change direction at a vertex in $T$, i.e. switch from following the direction of the edges to going against them and vice versa.
\medskip

The corresponding question for this model is to start a walk from $u$ following the direction of the edges and compare the probabilities for arriving at $v$ going with or against the direction of the last edge into $v$. For this model we may in fact prove the corresponding bunkbed theorem. Let $u_\rightarrow$ and $u_\leftarrow$ denote starting at $u$ following the direction of the edges (resp. going against the direction of the edges). If $u\in T$, then we can for both symbols start with or against the direction. Also let $\tillmed v$ and $\tillmot v$ denote entering $v$ going with (resp. against) the directions of the edge. Again, if $v\in T$ then it is in both cases legal to enter $v$ either going forward or  reverse direction.

\begin{theorem} \label{T:D2}
Let $G$ be any graph and $T\subseteq V(G)$. For any $u,v\in V(G)$ we have
\[ P(u_\rightarrow \hpil {D_2^{T}(G)} \tillmed{v})\ge P(u_\rightarrow \hpil {D_2^{T}(G)} \tillmot{v}).
\]
\end{theorem}

\begin{proof} First we consider all orientations of $G$ such that there is no walk from $u_\rightarrow$ to any vertex in $T$. In this case the right hand side is zero so the inequality is clear.

In the remaining cases we condition on the existence of a walk from $u_\rightarrow$ to some vertex in $T$. In this case we will construct a involution on the set of orientations which will show that the probability is equal arriving to $ \tillmed{v}$ and to $ \tillmot{v}$. 

To this end fix an orientation $O$ of $G$ and define $X(O)\subseteq V$ as all vertices $x$ to which there exists walks from $u_\rightarrow$ to both $ \tillmed{x}$ and $\tillmot{x}$. For instance every vertex on a directed path from $u_\rightarrow$ to a vertex in $T$ belongs to $X(O)$. This is because we may follow the path to the transversal vertex and then go backwards along the same path. Hence $X(O)\neq \emptyset$. If $v\in X(O)$ then we do nothing. If $v\notin X(O)$, let $F(O)\subseteq E$ be all edges between two vertices in $X(O)$. Now we define a new orientation $O^r$ by reversing the direction of all edges {\bf not} in $F(O)$. By construction $X(O)\subseteq X(O^r)$. If there were a vertex $x\in X(O^r)\backslash X(O)$, this would mean that there were two shortest paths $P_1, P_2$ in $G$ with orientation $O^r$ starting at some, possibly different, vertices in $X(O)$, using only edges in $E\backslash F(O)$ and ending in $\tillmed{x}$ and $\tillmot{x}$ respectively. But every edge on $P_1,P_2$ has the reverse orientation in $O$ than in $O^r$ and every vertex in $X(O)$ can be reached either way. Thus $P_1$ is a legal path also in orientation $O$ of $G$ but ending in $\tillmot{x}$ instead of $\tillmed{x}$, and the other way around for $P_2$. This gives a contradiction and we can conclude that $X(O)=X(O^r)$
and thus $F(O^r)=F(O)$ and $(O^r)^r=O$. 
There is a walk in $O$ from $u_\rightarrow$ to $ \tillmed{v}$ if and only if there is a walk from 
$u_\rightarrow$ to $ \tillmot{v}$ in $O^r$ and vice versa. The theorem follows.
\end{proof}

\medskip
\noindent\textbf{Model $D_3^{T}$:} Let $T\subseteq V$ . Every edge in $G$ is given one of the two possible directions with equal probability. A walk in $G$ may change direction at a vertex in $T$, i.e. switch from following the direction of the edges to going against them and vice versa. A walk must not use an edge in both directions.
\medskip

\begin{figure}[htb!] 
\center{
\epsfig{file=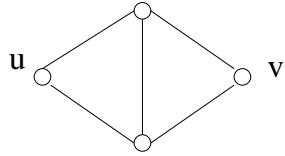, height=2cm}
\caption{A graph for which models $E_3^T$ and $D_3^T$ differ. Here $T=\{u,v\}$.}
\label{F:motexempel}
}
\end{figure}

The model $D_3^T$ seems closer to $E_3^T$ than $D_2^T$, but unfortunately they are not equivalent in general.
Figure \ref{F:motexempel} shows an example $G$ with four vertices and five edges, $T=\{u,v\}$, where 
$P(u_\rightarrow \hpil {D_3^{T}(G)} \tillmed{v})=P(u_\rightarrow \hpil {D_3^{T}(G)} \tillmot{v})=13/16$, whereas $P(u_0 \pil {E_3^{T}(G)} v_0)=P(u_0 \pil {E_3^{T}(G)} v_1)=7/8$.

We end with the corresponding bunkbed conjecture for model $D_3$.

\begin{conjecture}[$BBC_{D_3^T}$] 
Let $G$ be any graph and $T\subseteq V(G)$. For any $u,v\in V(G)$ we have
\[ P(u_\rightarrow \hpil {D_3^{T}(G)} \tillmed{v})\ge P(u_\rightarrow \hpil {D_3^{T}(G)} \tillmot{v}).
\]
\end{conjecture}

\section{A critical probability for finite graphs} \label{S:Critical}
We end this note with the definition of a critical probability for finite graphs that could be interesting to study further. Consider the following modification of Model $E_3$.

\medskip
\noindent\textbf{Model $E_5^{p,T}$:} Given a graph $G$ and $0\le p\le 1$, let $T\subseteq V(G)$ . Every edge in $G$ is colored red with probability $p$ and otherwise colored blue. A walk in $G$ may change color only at a vertex in $T$.
\medskip

Recall that we think of red edges as being downstairs (in the 0-layer) and blue as being upstairs.
Now we define the average probability that there is a walk from $u_0$ to $v_0$. That is, the walk must start from $u$ along a red edge (unless $u\in T$ then we can switch to a blue edge at once) and arrive to $v$ along a red edge (again unless $v\in T$).
\[P_G^p(u_0 \pil{} v_0):=\frac{1}{2^{|V|}}\sum_{T\subseteq V} P(u_0 \pil{E_5^{p,T}(G)} v_0).
\]
Similarly we define
\[P_G^p(u_0 \pil{} v_1):=\frac{1}{2^{|V|}}\sum_{T\subseteq V} P(u_0 \pil{E_5^{p,T}(G)} v_1).
\]
Intuitively it is clear that if $p$ is large (close to 1) the first quantity should be larger and vice versa if $p$ is close to 0. We conjecture that for any connected graph $G$ and any $u,v\in G$ there is a critical probability $p^c$ such that

\[P_G^{p}(u_0\pil{} v_0)
\left\{
\begin{array}{ccc}
 &< P_G^{p}(u_0\pil{} v_1), &\text{if $p<p^c$}    \\
 &= P_G^{p}(u_0\pil{} v_1), &\text{if $p=p^c$}    \\
 &> P_G^{p}(u_0\pil{} v_1), &\text{if $p>p^c$}    
\end{array}
\right.
\]
If this and the conjecture $BBC_3^{T}(G)$ are true, then $p^c<1/2$ for $G$. The inequality is strict because of the case $T=\emptyset$.

\medskip
\noindent
\textbf{Example:} Let $P_k$ be the path with $k$ edges and let $u,v$ be the endpoints. It is easy to compute that $p^c=1/3$  for $k=1$ and $p^c=\sqrt(11/12)-1/2$ for $k=2$. Defining an appropriate recursion one may also prove that the conjectured properties of $p^c$ holds for any path and that $p^c$ is increasing, monotone and converging to $1/2$ for  $k\longrightarrow \infty$. This may be interpreted as the endpoints of long paths being further apart. Does this make some sense also for other graphs?

\end{document}